\documentclass[11pt, reqno]{amsart}
\usepackage{amssymb,latexsym,amsmath,amsfonts,amsxtra}
\usepackage{mathrsfs}

\usepackage[all,cmtip]{xy}

\numberwithin{equation}{section}

\theoremstyle{definition}
\newtheorem{definition}{Definition}[section]

\theoremstyle{remark}
\newtheorem{remark}[definition]{Remark}

\theoremstyle{plain}
\newtheorem{theorem}[definition]{Theorem}
\newtheorem{result}[definition]{Result}

\newtheorem{proposition}[definition]{Proposition}

\newtheorem{corollary}[definition]{Corollary}

\newcommand{\nat}{\mathbb{N}}
\newcommand{\fib}{\boldsymbol{{\sf F}}}

\newcommand\partl[3]{\frac{\partial^{\raisebox{3pt}{$\scriptstyle {{#3}}$}}\!{#1}}{\partial{#2}}}

\newcommand{\OM}{\Omega}

\newcommand{\smoo}{\mathcal{C}}
\newcommand{\hol}{\mathcal{O}}

\newcommand\prj[1]{\pi_{\raisebox{-2pt}{$\scriptstyle {{#1}}$}}}

\newcommand{\bcdot}{\boldsymbol{\cdot}}

\newcommand{\lrarw}{\longrightarrow}

\newcommand{\cdim}{{\rm dim}_{\mathbb{C}}}
\newcommand{\varph}{\psi^y}

\newcommand{\cplx}{\mathbb{C}} 

\newcommand{\rea}{\mathbb{R}}

\begin{document}

\title[Rigidity of holomorphic maps]{Rigidity of holomorphic maps \\
between fiber spaces}

\author{Gautam Bharali}
\address{Department of Mathematics, Indian Institute of Science, Bangalore 560012, India}
\email{bharali@math.iisc.ernet.in}

\author{Indranil Biswas}
\address{School of Mathematics, Tata Institute of Fundamental Research, Mumbai 400005, India}
\email{indranil@math.tifr.res.in}

\thanks{The first author is supported in part by a UGC Centre of Advanced Study grant. The second author is supported
by a J.C. Bose Fellowship}



\begin{abstract}
In the study of holomorphic maps, the term ``rigidity'' refers to certain types of results
that give us very specific information about a general class of holomorphic maps owing to the
geometry of their domains or target spaces. Under this theme, we begin by studying when,
given two compact connected complex manifolds $X$ and $Y$, a degree-one holomorphic
map $f: Y\lrarw X$ is a biholomorphism. Given that the real manifolds underlying $X$ and
$Y$ are diffeomorphic, we provide a condition under which $f$ is a biholomorphism. Using
this result, we deduce a rigidity result for holomorphic self-maps of the total space of a holomorphic
fiber space. Lastly, we consider products $X=X_1\times X_2$ and $Y=Y_1\times Y_2$ of compact
connected complex manifolds. When $X_1$ is a Riemann surface of genus $\geq 2$, we 
show that any non-constant holomorphic map $F:Y\lrarw X$ is of a special form.
\end{abstract}
\maketitle

\section{Introduction}\label{S:intro}

The \textit{degree} of a continuous map $f:Y \longrightarrow X$ between
compact connected oriented smooth manifolds of dimension $n$ is defined as follows:
$$
f^*(1_X)\,=\, \text{degree}(f)\bcdot 1_Y\, ,
$$
where $1_X$ (respectively, $1_Y$) is the unique generator of $H^n(X,\,{\mathbb Z})$
(respectively, $H^n(Y,\,{\mathbb Z})$) compatible with the orientation. When $X$ and
$Y$ are compact connected complex manifolds and $f: Y\lrarw X$ is holomorphic, there
are natural situations in which, if $f$ is a degree-one map, then it is automatically
a biholomorphism. If $X=Y$ is a compact Riemann surface, this follows from
the Riemann--Hurwitz formula. We encounter a significant obstacle when
$\cdim{X} = \cdim{Y} \geq 2$. Note that $\text{degree}(f)$ as defined is precisely
the topological degree of of $f$. When $\text{degree}(f) = 1$, it follows that
the pre-image of a generic point in $X$ --- but \textit{not necessarily every} point --- is
a singleton. Thus, it may happen (and it does: consider the case when 
$f: Y \lrarw X$ is a blow-up) that there exists a non-empty proper subvariety of $Y$
on which $f$ fails to be injective. Hence, one must impose some conditions on $X$ and $Y$
for $f$ to be a biholomorphism. We explore this phenomenon when it is known {\em a priori}
that the {\em real} manifolds underlying $X$ and $Y$ are diffeomorphic. We prove the
following:

\begin{theorem}\label{T:degree1}
Let $X$ and $Y$ be compact connected complex manifolds such that the underlying
real manifolds are
diffeomorphic, and let $f : Y\lrarw X$ be a degree-one holomorphic map. If
$\dim H^1(X, \hol_X) = \dim H^1(Y, \hol_Y)$, then $f$ is a biholomorphism.
\end{theorem}

\begin{remark}
One expects some restrictions on the {\em complex} geometry of $X$ and $Y$ for a
degree-one map $f : Y\lrarw X$ to be biholomorphic. The assumption on cohomology in
Theorem~\ref{T:degree1} is a restriction of this sort. It is a rather mild condition:
it is satisfied, for instance, whenever $X$ and $Y$ are K{\"a}hler and real diffeomorphic.
Indeed, this follows immediately from the fact that for a compact K{\"a}hler manifold $Z$,
we have $\dim H^1(Z, \hol_Z)\,=\,\frac{1}{2}\dim H^2(Z,\,\cplx)$ (the latter is a consequence
of the Hodge decomposition). The vector space $H^1(X, \hol_X)$ parametrizes the space
of all infinitesimal deformations of any holomorphic line bundle on $X$ (the space of
holomorphic line bundles on $X$ is a group; infinitesimal deformations of these line bundles
are independent of the specific line bundle). So, the meaning of the cohomology condition in
Theorem~\ref{T:degree1} is that the infinitesimal deformations of any holomorphic
line bundle on $X$ and $Y$ are assumed to coincide.
This assumption is used in our proof essentially in this form.
\end{remark}

The above theorem forms the key final step in the following rigidity result for a holomorphic self-map
of a fiber space. Loosely speaking, if a holomorphic map of the total space sends just a single fiber
of a fiber space to another fiber, then it must be a map of fiber spaces. More precisely:

\begin{theorem}\label{T:fiberMap}
Let $X$ and $Y$ be complex manifolds, let $p:Y\lrarw X$ be a proper holomorphic surjective submersion
having connected fibers, and let $X$ be connected. Let $F: Y\lrarw Y$ be a holomorphic map such that
there exist points $a, b\in X$ with the property that $F$ maps the fiber $Y_a\,:=\,p^{-1}\{a\}$ into the
fiber $Y_b\,:=\,p^{-1}\{b\}$. Then:
\begin{enumerate}
 \item[$a)$] The map $F$ is a map of fiber spaces: i.e., there exists a holomorphic map 
 $f:X\lrarw X$ such that $p\circ F = f\circ p$.
 \item[$b)$] If, additionally, $F|_{Y_a}$ is a degree-one map from $Y_a$ to $Y_b$ and if
 $\dim H^1(Y_x, \hol_{Y_x})$ is independent of $x\in X$, then $F$ is a fiberwise biholomorphism.
\end{enumerate}
\end{theorem}

The term ``rigidity'' for holomorphic maps often refers to the phenomenon of a holomorphic
map being structurally simple owing to the geometry of its domain or target space; see, for
example, results by Remmert--Stein \cite[Satz\,12, 13]{remmertStein:eha60}, Kaup \cite[Satz\,5.2]{kaup:hkr68}
or Kobayashi \cite[Theorem\,7.6.11]{kobayashi:hcs98}. In the set-up of Theorem~\ref{T:fiberMap},
a rigidity result in this sense would require one to determine, for instance, conditions on $(Y,X,p)$
that would cause any $F:Y\lrarw Y$ to preserve at least one fiber.
This seems to be a difficult requirement. However, in the simpler set-up of certain product spaces,
we do get a rigidity result of the above-mentioned style. It comes as a corollary of the following:

\begin{proposition}\label{P:rs}
Let $Y=Y_1\times Y_2$ be a product of compact connected complex manifolds and let 
$f: Y\lrarw X$ be a holomorphic map into a compact Riemann surface of genus $\geq\,2$ . Then (denoting
each $y\in Y$ as $(y_1,y_2)$, $y_j\in Y_j, \ j=1, 2$) $f$ depends on at most one of $y_1$ and $y_2$.
\end{proposition}

This has the obvious corollary: 

\begin{corollary}\label{C:prod}
Let $Y=Y_1\times Y_2$ and $X=X_1\times X_2$ be products of compact connected complex manifolds.
Assume that  $X_1$ is a compact Riemann surface of genus $\geq\,2$ . Let $F=(F_1,F_2): Y\lrarw X$ be a 
holomorphic map. Then, there is a $j\in \{1,2\}$ such that $F$ has the form
$F(y_1, y_2)=(F_1(y_j), F_2(y_1,y_2))$.
\end{corollary}

One may compare the above to a result by Janardhanan
\cite{janardhanan:phmhpm11}. The hypothesis of Corollary~\ref{C:prod} is weakened to allow the factors
of $Y$ to have arbitrary dimension; but with $\cdim Y_j=\cdim X_j = 1$, Janardhanan is also able to handle
the case when some of the factors are non-compact.

\begin{remark}
The assumptions in Proposition~\ref{P:rs} and Corollary~\ref{C:prod} cannot be weakened
appreciably. Suppose $Y_1\neq Y_2$ are two holomorphically distinct compact
connected complex manifolds and suppose $X^\prime$ is a compact Riemann surface of
 genus $\geq\,2$ such that there are non-constant holomorphic maps
$f_j: Y_j\lrarw X^\prime, \ j=1, 2$. Then, just taking $X = X^\prime\times X^\prime$
and $f=(f_1,f_2)$ in Proposition~\ref{P:rs} shows that it cannot be true in general if
$\cdim{X}\geq 2$ (more involved examples can be constructed in which $\cdim{X}\geq 2$ and
is not a product). As for the requirement on the genus of $X$ being essential: the reader is referred
to \cite[Remark 1.7]{janardhanan:phmhpm11}.
\end{remark} 

\section{The proof of Theorem~\ref{T:degree1}}\label{S:degree1}

Let $n$ be the complex dimension of $X$ (also of $Y$). Let
\begin{equation}\label{fst}
f^*\, :\, H^*(X,\,{\mathbb Q})\, \longrightarrow\,
H^*(Y,\,{\mathbb Q})
\end{equation}
be the pullback homomorphism for $f$.
We will show that $f^*$ is an isomorphism. To do so,
given any non-zero class $c\in H^i(X,\,{\mathbb Q})$, take $c'\in H^{2n-i}(X,\,{\mathbb Q})$
such that $c\cup c'\not= 0$. Since
$$
((f^*c)\cup (f^*c'))\cap [Y]\,=\,f^*(c\cup c')\cap [Y]\,=\, 
\text{degree}(f)\bcdot (c\cup c')\cap [X]\, ,
$$
we conclude that $f^*c\,\not=\,0$ for $c\not= 0$.  Hence $f^*$ is injective.
We have $\dim H^*(X,\,{\mathbb Q})\,=\,\dim H^*(Y,\,{\mathbb Q})$
because $X$ and $Y$ are diffeomorphic
as real manifolds. Hence the injective homomorphism $f^*$ is an isomorphism.
\smallskip

The differential of $f$ produces a holomorphic section of the
holomorphic line bundle $\Omega^n_Y\otimes f^*(\Omega^n_X)^*$ on $Y$. This
section will be denoted by $s$. Let
\begin{equation}\label{fst2}
D\,:=\,\text{Divisor}(s)\,\subset\,Y
\end{equation}
be the effective divisor of this section. We note that $f$ is a biholomorphism from
$Y\setminus D$ to $X\setminus f(D)$. So to prove the theorem, it suffices
to show that $D$ is the zero divisor.
\smallskip

To prove that $D$ is the zero divisor, we first note that $f(D)\subset X$ is
of complex codimension at least two. Indeed, the given condition that the degree
of $f$ is one implies that if, for an irreducible component $D'$ of $D$, the image $f(D')$
is a divisor in $X$, then $f$ is an isomorphism on a neighborhood of $D'$. But this
implies that $D'$ is not contained in $D$. Therefore,
$f(D)\subset X$ is of complex codimension at least two.
\smallskip

Let $c(D)\in H^2(Y,\,{\mathbb Q})$ be the cohomology class of $D$.
Since $f(D)\subset X$ is of complex codimension at least two, and $f^*$ in
\eqref{fst} is an isomorphism, it follows that
\begin{equation}\label{cd0}
c(D)\,=\,0\, .
\end{equation}

Now consider the short exact sequence of sheaves on $Y$
$$
0\,\longrightarrow\, {\mathbb Z}\,\longrightarrow\, {\mathcal O}_Y\,\stackrel{\exp}
{\longrightarrow}\, {\mathcal O}^*_Y\,\longrightarrow\, 0\, .
$$
Let
\begin{equation}\label{e3}
H^1(Y,\,{\mathcal O}_Y)\,\stackrel{\beta}{\longrightarrow}\, H^1(Y,\,{\mathcal O}^*_Y)\,\stackrel{q}
{\longrightarrow}\, H^2(Y,\,{\mathbb Z})
\end{equation}
be the long exact sequence of cohomologies associated to it. The class in
$H^1(Y, {\mathcal O}^*_Y)$ for the holomorphic line bundle $\Omega^n_Y\otimes
f^*(\Omega^n_X)^*$ (this is the line bundle associated to $D$) will be denoted by $\gamma(D)$. For
any positive integer $m$, the class in $H^1(Y, {\mathcal O}^*_Y)$ for the holomorphic line
bundle $(\Omega^n_Y\otimes f^*(\Omega^n_X)^*)^{\otimes m}$ is $\gamma(D)^m$. From
\eqref{cd0} if follows that there is a positive integer $N$ such that
$$
q(\gamma(D)^N)\,=\,0\, ,
$$
where $q$ is the homomorphism in \eqref{e3}. Therefore, there is a cohomology class
$\alpha\in H^1(Y, {\mathcal O}_Y)$ such that
\begin{equation}\label{clv}
\beta(\alpha)\,=\, \gamma(D)^N.
\end{equation}

Consider the pullback homomorphism
\begin{equation}\label{d-F}
F\,:\,H^1(X,\,{\mathcal O}_X)\,\longrightarrow\,
H^1(Y,\,f^*{\mathcal O}_X)\,=\,H^1(Y,\,{\mathcal O}_Y)
\end{equation}
for $f$. We will show that $F$ is injective. To prove this, first note that
$f_*{\mathcal O}_Y = {\mathcal O}_X$. Using this isomorphism,
the natural homomorphism
\begin{equation}\label{d-F2}
H^1(X,\,f_*{\mathcal O}_Y)\,\longrightarrow\, H^1(Y,\,{\mathcal O}_Y)
\end{equation}
coincides with $F$ in \eqref{d-F}. But the homomorphism in \eqref{d-F2}
is injective. Hence $F$ is injective. We now invoke the assumption that
$\dim H^1(X, {\mathcal O}_X)\,=\,\dim H^1(Y, {\mathcal O}_Y)$.
Since $F$ is an injective homomorphism between vector spaces of same
dimension, we conclude that $F$ is an isomorphism.
\smallskip

Take $\alpha\in H^1(X, {\mathcal O}_X)$ such that
\begin{equation}\label{Fp}
\alpha\,=\,F(\alpha')\, ,
\end{equation}
where $\alpha$ is the cohomology class in \eqref{clv}. Let $L$ be the
holomorphic line bundle on $X$ corresponding to the element
$\beta'(\alpha')\in H^1(X, {\mathcal O}^*_X)$, where
$$
\beta'\,:\,H^1(X,\,{\mathcal O}_X)\,\longrightarrow\,
H^1(X,\,{\mathcal O}^*_X)
$$
is the homomorphism given by $\exp : {\mathcal O}_X \longrightarrow
{\mathcal O}^*_X$. From
\eqref{clv} and \eqref{Fp} it follows that
$$
(\Omega^n_Y\otimes f^*(\Omega^n_X)^*)^{\otimes N}\,=\,f^*L\, .
$$

The holomorphic line bundle
$(\Omega^n_Y\otimes f^*(\Omega^n_X)^*)^{\otimes N}$ is holomorphically trivial
on $Y\setminus D$. Consequently, $L$ is holomorphically trivial on $X\setminus f(D)$
(recall that $f$ is a biholomorphism from $Y\setminus D$ to $X\setminus f(D)$).

Since $L$ is holomorphically trivial on $X\setminus f(D)$, and $f(D)\subset X$
is of complex codimension at least two, the holomorphic line bundle $L$ on $X$ is
holomorphically trivial. Therefore, the holomorphic line bundle $(\Omega^n_Y
\otimes f^*(\Omega^n_X)^*)^{\otimes N}\,=\,f^*L$ on $Y$ is holomorphically trivial.
\smallskip

Finally, consider the section $s$ in \eqref{fst2}. Note that $s^{\otimes N}$ is a holomorphic section
of $(\Omega^n_Y\otimes f^*(\Omega^n_X)^*)^{\otimes N}$ vanishing on
$D$ and nonzero elsewhere. On the other hand, any holomorphic section of the
holomorphically trivial line bundle on $Y$ is either nowhere zero, or it is identically
zero. Therefore, we conclude that $D$ is the zero divisor. This completes the proof
of Theorem~\ref{T:degree1}. \qed

\section{The proof of Theorem~\ref{T:fiberMap}}\label{S:fiberMap}

Given the conditions on $p$, it follows from Ehresmann's theorem \cite{ehresmann:ciefd50} that the
triple $(Y,X,p)$ is a $\smoo^\infty$-smooth fiber bundle with fiber $\fib$. Thus, for each $x\in X$, there is
a connected open neighborhood $U^x$ of $x$ and a diffeomorphism
$\varphi_x:p^{-1}(U^x)\lrarw U^x\times\fib$ such that the diagram

\begin{center}
$\;$
 \xymatrix{
 p^{-1}(U^x)	\ar[r]^{\varphi_x}
				\ar[rd]_{p}
 &U^x\times\fib	\ar[d]^{{\sf proj}_1} \\
 &U^x}
\end{center}

\noindent{commutes (here, ${\sf proj}_1$ denotes the projection of
$U^x\times\fib$ onto the first factor).}
\smallskip

Write $S:=\{x\in X: \exists\, x^\prime\in X \ \text{such that} \ F(Y_x)\subset Y_{x^\prime}\}$. We shall
first show that $S$ is an open set. Suppose $x_0\in S$, whence there is a point $x^\prime\in X$ such that
$F(Y_{x_0})\subset Y_{x^\prime}$. Let $B_{x^\prime}$ be a neighborhood around $x^\prime$
that is biholomorphic to a ball. By continuity of $F$ and $p$, and due to compactness of $Y_{x_0}$, we
can find an open neighborhood $V_{x_0}$ of $x_0$ such that 
$p\circ F(p^{-1}(V_{x_0}))\subset B_{x^\prime}$. Recall that, by hypothesis, the fibers of
$p$ are connected. Thus, for each $x\in V_{x_0}$:
\begin{itemize}
 \item $Y_x$ is a connected, compact complex manifold.
 \item $p\circ F$ maps each $Y_x$ into an Euclidean open set.
\end{itemize}
It follows from the maximum modulus theorem that $\left. p\circ F\right|_{Y_x}$ is constant for each 
$x\in V_{x_0}$. This means: $x_0\in S \ \Rightarrow \ V_{x_0}\subset S$. In other words, 
$S$ is an open set. 
\smallskip

We now argue that $S$ is closed. Write $d=\cdim X$ and $k=\cdim Y$. There is a 
$\smoo^\infty$-atlas $\{(W^y;x^y_1,\dots,x^y_{2d},\Phi^y_1,\dots,\Phi^y_{2k}):y\in Y\}$ of
$Y$, where $W^y$ is a coordinate patch centered at $y$ such that
\[
 Y_{p(z)}\cap W^y \ = \ \{x^y_1=C_1^{p(z)},\dots,x^y_{2d}=C_{2d}^{p(z)}\} \; \; \;
 \forall\, z\in W^y,
\]
and where each $C_j^{p(z)}\in \rea$ is a constant that depends only on $p(z)$. More
descriptively: $(x^y_1,\dots,x^y_{2d},\Phi^y_1,\dots,\Phi^y_{2k})$ imposes a
smooth product structure upon $W^y$ in such a way that each $\rea$-affine
slice of the domain $(x^y_1,\dots,x^y_{2d},\Phi^y_1,\dots,\Phi^y_{2k})(W^y) =: G^y\subset
\rea^{2(d+k)}$ by a translate of the ``$\Phi^y_1\dots\Phi^y_{2k}$-plane'' parametrizes 
a patch of some fiber of $p$. Let
\[
 \OM_y \ := \ \text{the connected component of $W^y\cap F^{-1}(W^{F(y)})$ containing $y$.}
\]
We claim that
\begin{equation}\label{E:closedLoc}
 p^{-1}(S)\cap \OM_y \ = \ \bigcap_{i=1}^{2d}\,\bigcap_{\alpha\in \nat^{2k}\setminus\{(0,\dots,0)\}}
 \negthickspace\!\bigg\{z \in \OM_y: \partl{\big(x^{F(y)}_i\circ F\big)}{\Phi^{y,\alpha}}{\alpha}(z)=0\bigg\}.
\end{equation}
Our notation $\partial^\alpha g/\partial\Phi^{y,\alpha}$ --- for any  $\smoo^\infty$-smooth
function $g: \OM_y\lrarw \rea$ --- perhaps needs a little
clarification. Denote by $\varph$ the coordinate map
$\varph := (x^y_1,\dots,x^y_{2d},\Phi^y_1,\dots,\Phi^y_{2k}): \OM_y\lrarw G^y$ described
above. Let us write $\varph(z) = (x_1,\dots,x_{2d},\Phi_1,\dots,\Phi_{2k})$, which are just
points varying through $G^y$. Then, $\partial^\alpha g/\partial\Phi^{y,\alpha}$ is
defined as: 
\[
 \partl{g}{\Phi^{y,\alpha}}{\alpha}(z)\,:=\,\frac{\partial^{|\alpha|}{g\circ(\varph)^{-1}}}
 {\partial\Phi_1^{\alpha_1}\dots\partial\Phi_{2k}^{\alpha_{2k}}}(\varph(z)), \; \; \; z\in \OM_y.
\]
That
$p^{-1}(S)\cap \OM_y$ is a subset of the set on the right-hand side of \eqref{E:closedLoc} (call it
$K_y$) is clear. Now, if $z\in K_y$, it implies that there is a small open neighborhood $N_z$ of $z$
such that the holomorphic map 
$\left.p\circ F\right|_{Y_{p(z)}}:Y_{p(z)}\lrarw X$ is constant on the set $N_z\cap Y_{p(z)}$, since the 
Taylor coefficient (relative to local coordinates) of $\left.p\circ F\right|_{Y_{p(z)}}$ at $z$ corresponding
to each $\alpha\neq (0,\dots,0)$ is zero. By the principle of
analytic continuation, $p\circ F$ must be constant on $Y_{p(z)}$. So, $p(z)\in S$ and hence
$K_y\subset p^{-1}(S)\cap \OM_y$. This establishes \eqref{E:closedLoc}. 
\smallskip

By \eqref{E:closedLoc}, the intersection $p^{-1}(S)\cap \OM_y$ is closed relative to each $\OM_y$. As 
$\{\OM_y:y\in Y\}$ is an open cover of $Y$, we see that $p^{-1}(S)$ is closed. It is now easy
to see, as $(Y,X,p)$ is locally trivial, that $S$ is closed. By hypothesis, $S\neq \varnothing$. As
$X$ is connected, it follows that $S=X$.
\smallskip

Since $S=X$, the map $f:X\lrarw X$ given by
\[
 f(x) \ := \ \text{$p\circ F(y)$ for any $y\in Y_x$}
\]
is well-defined. Clearly $p\circ F = f\circ p$.
Let us now {\em fix} a point $p\in \fib$. For any $x_0\in X$,
\[
 f|_{U^{x_0}}(x) \ = \ p\circ F\circ\varphi_{x_0}^{-1}(x,p) \; \; \; \forall x\in U^{x_0},
\]
whence $f$ is $\smoo^\infty$-smooth. Therefore, relative to any local holomorphic coordinate system,
we can apply the Cauchy--Riemann operator to the equation $p\circ F = f\circ p$ to conclude that $f$ is
holomorphic. This establishes part\,$(a)$ of Theorem~\ref{T:fiberMap}.
\smallskip

To prove part $(b)$, recall that given any set $\Sigma\subset \smoo(\fib;\fib)$ --- where the latter denotes
the space of all continuous self-maps of $\fib$ endowed with the compact-open
topology --- the function from
$\Sigma$ to ${\mathbb Z}$ defined by
$\psi\longmapsto {\rm degree}(\psi)$ is locally constant. Hence, as
$X$ is connected and, by hypothesis, ${\rm degree}(F|_{Y_a})=1$, we have
\begin{equation}\label{E:degCond}
 {\rm degree}(F|_{Y_x}) \ = \ 1 \; \; \; \forall\, x\in X.
\end{equation}
From our hypothesis, it follows that
\[
 \dim H^1(Y_x, \hol_{Y_x}) \ = \ \dim H^1(Y_{f(x)}, \hol_{Y_{f(x)}})
\; \; \; \forall\, x\in X.
\]
In view of this, \eqref{E:degCond}, and the fact that fibers are connected, we may apply 
Theorem~\ref{T:degree1} to conclude that $F$ is a fiberwise biholomorphism. This
completes the proof of Theorem~\ref{T:fiberMap}. \qed

\section{Concerning Proposition~\ref{P:rs} and Corollary~\ref{C:prod}}\label{S:prod}

Corollary~\ref{C:prod} is an absolutely obvious consequence of Proposition~\ref{P:rs}. Thus,
we shall only discuss the latter.
\smallskip

The proof of Proposition~\ref{P:rs} relies upon the following result of Kobayashi--Ochiai:

\begin{result}[Theorem\,7.6.1, \cite{kobayashi:hcs98}]\label{R:KobO}
Let $X$ and $Y$ be two compact complex-analytic spaces. If $Y$ is of general type, then
the set of all dominant meromorphic maps of $X$ to $Y$ is finite.
\end{result}

A couple of remarks are in order. First: we shall not define the term {\em general type}
here as it is somewhat involved. We refer the reader to \cite[\S\,7.4]{kobayashi:hcs98}. The
fact that we need from \cite[\S\,7.4]{kobayashi:hcs98} is that any compact Riemann surface with genus
$\geq\,2$ is of general type.
\smallskip

Secondly: in their original announcement and proof of the above
result in \cite{kobayashiOchiai:mmccsgt75}, Kobayashi and Ochiai require $X$ to be
Moi{\v{s}}ezon. However, in a footnote to \cite{kobayashiOchiai:mmccsgt75}, the authors
observe that this restriction on $X$ can be removed. The relevant argument is
presented in the proof of \cite[Theorem\,7.6.1]{kobayashi:hcs98}.
\medskip 

\noindent{\em The proof of Proposition~\ref{P:rs} .} 
Since there is nothing to prove if $f$ is constant, we shall assume that $f$ is non-constant.
For each $y=(y_1,y_2)\in Y$, we define:
\[
 f^{y_2} \ := \ f(\,\bcdot\,,y_2),  \qquad	f^{y_1}	\ :=  \ f(y_1,\,\bcdot\,).
\]
Assume that both $f^{y_2}$ and $f^{y_1}$ are  non-surjective for each $y\in Y$.
As $f^{y_2}$ is a holomorphic map between {\em compact} complex manifolds, for each fixed $y_2\in Y_2$, 
$f^{y_2}(Y_1)\varsubsetneq X$ is a complex-analytic subvariety. Thus, $f^{y_2}$ is a constant map
for each $y_2\in Y_2$; call this constant $c(y_2)$. Fix a point $a\in Y_1$. By the same argument as
above, we have
\[
 c(y_2) \ = \  f(a, y_2) \ = \ C \quad \forall y_2\in Y_2,
\]
where $C$ is a constant. This contradicts the fact that $f$ is non-constant. Hence, there exists a 
$j\in \{1,2\}$ and a point $a=(a_1, a_2)\in Y$ such that $f^{a_j}$ is surjective.
\smallskip
 
There is no loss of generality in taking $j=1$. By continuity of $f$ and compactness of
$Y_2$, it follows that there is an open neighborhood $U\subset Y_1$ of $a_1$ such that the maps
$f(z,\,\bcdot\,) : Y_2\lrarw X$ are surjective for each $z\in U$. Since $X$ has genus $\geq 2$, it is of
general type. It follows from Result~\ref{R:KobO} (which is a refinement of an argument of
Kobayashi--Ochiai \cite{kobayashiOchiai:mmccsgt75}) that the surjective holomorphic maps from $Y_2$
to $X$ form a  finite set. Hence, $f$ restricted to the open set $\prj{Y,\,1}^{-1}(U)$ (where $\prj{Y,\,k}$ denotes
the projection of $Y$ onto its $k$th factor) is independent of $y_1$. By the principle of analytic continuation,
it follows that $f$ is independent of $y_1$. Hence the desired result. \qed
\medskip

\noindent{
\textbf{Acknowledgements.}\,We are very grateful to N. Fakhruddin for his help
with the proof of Theorem~\ref{T:degree1}. We also thank the referee of this work
for helpful comments on our exposition.}

\end{document}